\documentclass[11pt]{article}

\usepackage{amsmath}
\usepackage{amssymb} 
\usepackage{amscd} 
\usepackage{theorem} 
\def\oplusinf{\mathop{\oplus}} 
\def\id{{\mbox{Id}}}
\def\im{{\mbox{Im}}}
\def\dim{{\mbox{dim}}}
\def\ker{{\mbox{Ker}}}

\def\cala{{\cal A}} 
 
\def\call{{\cal L}}

\def\calh{{\cal H}}

\def\calc{{\cal C}} 
\def\cale{{\cal E}}
\def\calu{{\cal U}} 
\def\calf{{\cal F}}
\def\calr{{\cal R}}
\def\fraca{{\frak A}}
 \def\fraci{{\frak I}}
 \def\fracg{{\frak g}} 
\def\bbbone{\mbox{\rm 1\hspace {-.6em} l}}
\def\ue{\underline{E}}
\def\uf{\underline{F}}
\def\uk{\underline{K}}
\def\ux{\underline{X}}

\def\Section#1{\section{#1}\setcounter{equation}{0}}

\newtheorem{theorem}{THEOREM}
\newtheorem{lemma}{LEMMA} 
\newtheorem{proposition}{PROPOSITION}

\begin{document}
\baselineskip=0,75cm
\begin{center} 
{\Large\bf GENERALIZED HOMOLOGIES}

\vspace{3mm}

{\Large\bf FOR THE ZERO
MODES}

\vspace{3mm}

{\Large\bf  OF THE $SU(2)$ WZNW MODEL} 

\end{center} 
\vspace{0.75cm}

\begin{center} Michel DUBOIS-VIOLETTE \\
\vspace{0.3cm} {\small Laboratoire de Physique Th\'eorique
\footnote{Unit\'e Mixte de Recherche du Centre National de la
Recherche Scientifique - UMR 8627}\\ Universit\'e Paris XI,
B\^atiment 210\\ F-91 405 Orsay Cedex, France\\
patricia$@$lyre.th.u-psud.fr}

\vspace{0,5cm}

and

\vspace{0,5cm}

Ivan T. TODOROV\\
\vspace{0.3cm} {\small 
Institute for Nuclear Research and Nuclear Energy\\ Bulgarian
Academy of Sciences\\ BG-1784 Sofia, Bulgaria
\\

todorov$@$inrne.bas.bg}

\end{center} \vspace{1cm}

\begin{center} \today \end{center}

\vspace {1cm}

\begin{center}

{\sl Dedicated to the memory of Mosh\'e Flato.}
\end{center}

\vspace {3cm}

\noindent L.P.T.-ORSAY 99-13
\newpage

\begin{abstract}
We generalize the BRS method for the (finite-dimensional) quantum
gauge theory involved in the zero modes of the monodromy
extended $SU(2)$ WZNW model. The generalization consists of a
nilpotent operator $Q$ such that $Q^h=0$ ($h=k+2=2,3,\dots$ being
the height of the current algebra representation) acting on an
extended state space. The physical subquotient is identified with the direct
sum $\displaystyle{\oplusinf^{h-1}_{n=1}}\ker(Q^n)/\im(Q^{h-n})$.
\end{abstract}

\Section{Introduction}

The Wess-Zumino-Novikov-Witten (WZNW)model was originally
formulated  \cite{EdW} in terms of a multivalued action for a
field $g(x,t)$ (which maps  the cylinder $S^1\times \mathbb R$  
into a
Lie group $G$). Its solution \cite{KZ} came  however within the
axiomatic (or ``bootstrap") approach making use of  the
representation theory of affine Kac-Moody algebras. This
solution  exhibits some puzzling features like the appearance
of non-integer  (``quantum") statistical dimensions (which
appear as positive real  solutions of the fusion rules'
equations) while the corresponding 2-dimensional (2D) fields satisfy
local ``Bose type" commutation relations. The gradual
understanding of these features only began with the 
development of the canonical approach to the model (\cite{F}, 
 \cite{G}, \cite{AF}, \cite{FG}-\cite{FHT3})
in which one splits the field $g$ into left and
right movers' (chiral) components:
 
  \begin{equation}
g^A_B(x,t)=u^A_\alpha(x-t) \bar u^\alpha_B(x+t)=u^A_\alpha(x^-)
\bar u^\alpha_B(x^+)
\label{1.1}
\end{equation}
These chiral components reveal a hidden quantum group symmetry 
(under ``gauge
transformations" $u^A_\alpha\mapsto u^A_\sigma
T^\sigma_\alpha$ with non-commuting entries). The phase space of
the  theory is extended (compared to the 2D gauge invariant
construction)  by the chiral zero modes, including the monodromy
degrees of freedom  that appear in the twisted periodicity
condition (corresponding to $g(x+2\pi,t)=g(x,t)$) 
  \begin{equation}
u(x+2\pi)=u(x)M, \ \ \bar u(x+2\pi)=\bar M^{-1}\bar u(x).
\label{1.3}
\end{equation}
The resulting extended WZNW model is understood on
the classical  level \cite{F},\cite{G} while its quantization has only
been attempted in a  lattice approach \cite{FG} and has not been
brought to a form yielding a  satisfactory continuum limit. The
direct investigation of the  quantum model
\cite{FHT2},\cite{FHT3} has singled out a nontrivial
finite-dimensional  gauge theory problem. The present paper
constructs the physical state-space of the zero modes for
$G=SU(2)$ in terms of generalized Becchi-Rouet-Stora
\cite{BRS} (BRS) generalized homologies thus providing a complete solution to this 
problem.\\
We proceed to describing the problem in more detail
and to outlining  the content of the paper. 
Following
\cite{FHT2}, \cite{FHT3}, let us expand $u^A_\alpha$ and $\bar
u^\alpha_B$ into {\sl chiral vertex operators} $u^A_i$, $\bar
u^i_B$ which diagonalize the monodromy:
  \begin{equation}
u^A_\alpha(x)=a^i_\alpha u^A_i(x,p),\ \ \ \bar
u^\alpha_B(y)=\bar a^\alpha_j \bar u^j_B(y,\bar p)
\label{1.5}
\end{equation}
(the repeated indices $i,j$ are summed from  1 to 2,
$a^{1,2}$ standing for $a^{+,-}$ of \cite{D-VT} - cf.
\cite{HPT}).
Here $p$ is the shifted Lie algebra weight ($p=\Lambda+\rho$)
which labels the monodromy eigenvalues\footnote{Such a simple picture only works for integrable
highest weights, $0<p<h$. Going beyond this limit requires
dealing with indecomposable representations of the $su(2)$
current algebra involving singular vectors in the associated
Verma modules. The monodromy matrix would then contain a Jordan
cell for each pair $(p,2h-p;\  0<p<h)$ of highest weights.}.
The irreducible representations of the  quantum
universal enveloping algebra $U_q(\mathfrak{sl}_2)$ are labelled by
$q^{\pm p}$ with  $q=e^{i\frac{\pi}{h}}$,where $h=2,3,\dots$ 
is the {\sl height} of the associated current algebra representation 
($h=k+2$ where $k$ is the Kac-Moody level).The nontrivial
finite-dimensional problem singled out by
the canonical (hamiltonian) approach involves a
pair of {\sl quantum matrix algebras} $\fraca_h\otimes \bar\fraca_h$
generated by the {\sl zero mode vertex operators} $a^i_\alpha$
and $\bar a^\alpha_j$ (see \cite{HPT}) which reflects the
basic properties of composition and braiding of current algebra
modules\footnote{A BRS treatment of the Wakimoto module
corresponding to the factor  $u^A_i(x,p)$ in (1.3) is contained in
\cite{BF}.}.
The problem is to develop a ``$q$-gauge theory"
approach that would allow to extract an ($h-1$)-dimensional 
(generalized) BRS cohomology from the $h^4$-dimensional
$q$-Fock space module $\calh=\calf\otimes \bar\calf$ of
$\fraca_h\otimes \bar\fraca_h$.
A step towards its solution was made in \cite{D-VT}. After
singling out a $(2h-1)$-dimensional subspace $\calh_I$ of quantum
group invariant vectors of $\calh$, we proved  that the 
nilpotent operator
  \begin{equation}
A=a^2_\alpha \bar a^\alpha_2 
\label{1.7}
\end{equation}
 satisfying $A^h=0$ in $\calh$ and $A(\calh_I)\subset \calh_I$
has one-dimensional generalized homologies
$H_{(n)}(\calh_I,A)=\ker(A^n:\calh_I\rightarrow
\calh_I)/A^{h-1}(\calh_I)$ on $\calh_I$, $n\in \{1,\dots,h-1\}$.
The direct sum
$\displaystyle{\oplusinf^{h-1}_{n=1}}H_{(n)}(\calh_I,A)$
was then identified with the $(h-1)$-dimensional physical 
subquotient.\\
The objective of this paper is to extend this construction in
such a way that quantum group invariance has not to be imposed
as an extra constraint (in other words,
we solve the first problem stated in the concluding Section 3 of
\cite{D-VT}). It turns out that to this end it is {\sl necessary}
(as demonstrated in Section 3) to extend (in a suitable way) the
space $\calh=\calf\otimes \bar \calf$; this fact is not
completely unexpected since it corresponds in the usual
situation to the addition of ghost's states.
A {\sl minimal} canonical
construction achieving our goal is presented in Section 4 and is
related in Section 5 to a (generalized) Hochshild complex (see
\cite{D-V}-\cite{D-VK}). It is worth noticing here
that the construction relies on a generalization, in the context
of generalized (co)-homology, of an elementary spectral
sequence's argument. Finally, in Section 6 we identify the
generalized homology of $(\calh_I,A)$ as a part of a generalized
homology of Hochschild cochains and we compare our
constructions and results with the BRS-like ones.\\
Throughout the paper $h$ is an integer greater or equal to 2
and $q=\exp(i\frac{\pi}{h})$. By an $h$-{\sl differential vector
space}, we  mean a vector space $E$ equipped with a nilpotent
endomorphism $d$, its $h$-{\sl differential}, satisfying
$d^h=0$. The generalized homology of $(E,d)$ is then the family
of vector spaces $H_{(k)}(E,d)=\ker(d^k)/\im(d^{h-k})$, $k\in
\{1,\dots,h-1\}$. An $h$-{\sl complex} will be here an
$h$-differential vector space which is $\mathbb Z$-graded and
such that its differential is of degree 1. If $(E,d)$ is an
$h$-complex with $E=\displaystyle{\oplusinf_n}E^n$, then the
$H_{(k)}(E,d)$ are graded
$H_{(k)}(E,d)=\displaystyle{\oplusinf_n}H^n_{(k)}(E,d)$ with
$H^n_{(k)}(E,d)=\ker(d^k:E^n\rightarrow
E^{n+k})/d^{h-k}(E^{n+k-h})$. More generally we use the notations
of \cite{D-V2} for generalized (co)homology. Concerning
$q$-numbers, we use here the convention of \cite{D-VT}, that is 
$[n]=\frac{q^n-q^{-n}}{q-q^{-1}}$
which differs in several respects from the one of \cite{D-V2}.
In the next section we give a summary
of relevant parts of earlier work (\cite{FHT3}, \cite{D-VT},
\cite{HPT}).

\Section{Background, preliminaries}

The quantum matrix algebra $\fraca$ is characterized by
$R$-matrix exchange relations and a determinant condition
\cite{HPT}. The exchange relations are written conveniently in
terms of a pair of quantum antisymmetrizers $\cala$ and
$\cala(p)$:
\begin{equation}
a_1a_2\cala=\cala(p)a_1a_2,\ \ \mbox{i.e.}\ \ 
a^{i_1}_{\sigma_1}a^{i_2}_{\sigma_2}
\cala^{\sigma_1\sigma_2}_{\alpha_1\alpha_2}=\cala(p)^{i_1i_2}_{s_1s_2}
a^{s_1s_2}_{\alpha_1\alpha_2};
\label{1.9}
\end{equation}
here
\[
\begin{array}{l}
\cala^{\alpha_1\alpha_2}_{\beta_1\beta_2}=
\cale^{\alpha_1\alpha_2}\cale_{\beta_1\beta_2}=
q^{\varepsilon_{\alpha_2\alpha_1}}
\delta^{\alpha_1\alpha_2}_{\beta_1\beta_2}
-\delta^{\alpha_1\alpha_2}_{\beta_2\beta_1}
\left(\delta^{\alpha_1\alpha_2}_{\beta_1\beta_2}\equiv
\delta^{\alpha_1}_{\beta_1}\delta^{\alpha_2}_{\beta_2}\right)
\\
\\
\varepsilon_{\alpha\beta}=\left\{\begin{array}{rll}
1 & \mbox{for} & \alpha>\beta\\
0 & \mbox{for} & \alpha=\beta\\
-1 & \mbox{for} & \alpha<\beta
\end{array}\right.,\ \ 
\left(\cale^{\alpha\beta}\right)=\left(
\begin{array}{cc}
0 & -q^{1/2}\\
q^{-1/2} & 0
\end{array}
\right)= 
\left(\cale_{\alpha\beta}\right),\\
\\
\cala(p)^{i_1i_2}_{j_1j_2}=
\frac{[p+i_1-i_2]}{[p]}\varepsilon^{i_1i_2}
\varepsilon_{j_1j_2}=\frac{[p+i_1-i_2]}{[p]}
\left(\delta^{i_1i_2}_{j_1j_2}-\delta^{i_1i_2}_{j_2j_1}\right);
\end{array}
\]
both $\cala$ and $\cala(p)$ satisfy the Hecke algebra condition
\[
\cala^2=[2]\cala, \ \ \cala(p)^2=[2]\cala(p)\ \
([2]=q+q^{-1})
\]
and a braiding (Temperley-Lieb-Martin) property 
\[
\cala_{12}\cala_{23}\cala_{12}-\cala_{12}=
\cala_{23}\cala_{12}\cala_{23}-\cala_{23}=0
\left((\cala_{12})^{\alpha_1\alpha_2\alpha_3}_{\alpha_1\alpha_2\alpha_3}
=
\cala^{\alpha_1\alpha_2}_{\beta_1\beta_2}\delta^{\alpha_3}_{\beta_3}
\dots\right)
\]
which involves a change of $p$ for $\cala_{23}(p)$ 
(see\cite{HPT}); 
in particular,
\[
a^i_2a^i_1=q a^i_1a^i_2, \ \ \ [a^1_\alpha,a^2_\alpha]=0.
\]
The algebra $\fraca$ contains, by definition, the field $\calr$
of rational functions of $q^p$ (which enter the exchange
relations (\ref{1.9})), $a^i_\alpha$ shifting $p$ according to the
law
\[
q^pa^1_\alpha=a^1_\alpha q^{p+1},\ \ q^pa^2_\alpha=a^2_\alpha
q^{p-1}.
\]
The determinant condition allows to express a quadratic
combination of $a^i_\alpha$ as a function of $q^p$
\[
\det a:=\frac{1}{[2]} \varepsilon_{ij}a^i_\alpha
a^j_\beta\cale^{\alpha\beta}=[p]\Rightarrow \varepsilon_{ij}a^i_\alpha
a^j_\beta=[p]\cale_{\alpha\beta}.
\]
Identical relations are satisfied by $\bar a^\alpha_j$; we have,
in particular,
\[
\bar a^2_j\bar a^1_j=q\bar a^1_j\bar a^2_j,\ [\bar a^\alpha_1,
\bar a^\alpha_2]=0,\ q^{\bar p}\bar a^\alpha_{1,2}=
\bar a^\alpha_{1,2}q^{\bar p \pm 1},\varepsilon^{ij}\bar a^\alpha_i
\bar a^\beta_j=[\bar p]\cale^{\alpha\beta}.
\]
The algebra $\fraca$ admits a (two-sided) ideal
$\fraci_h$  generated by $(a^i_\alpha)^h\ \alpha,i=1,2$
and $[hp]$.
The factor algebra 
$\fraca_h=\fraca/\fraci_h$
is finite dimensional. It admits an $h^2$-dimensional Fock
space module $\calf$ with basis
\[
\begin{array}{l}
\vert p,m>=(a^1_1)^m(a^1_2)^{p-1-m}\vert1,0>,\ \ \mbox{where}\ \
a^2_\alpha \vert 1,0>=0,\\
\\
0<p<2h\ \ \mbox{and}\ \  \  \mbox{max}\ (0,p-h)\leq m\leq \ 
\mbox{min}\ (p-1,h-1).\ \
\end{array}
\]
Similar statements are valid for $\bar\fraca$.\\
The $h^4$ dimensional tensor product space $\calh=\calf\otimes
\bar \calf$ carries a representation of a tensor product of
quantum universal enveloping algebras (QUEA) which we proceed to
define.\\
The first factor is the diagonal $U_q(\mathfrak{sl}_2)$ related to the
left and right monodromies $M$ and $\bar M$ (appearing in
(\ref{1.3})) as follows. Denoting by $\Delta$ the $U_q(\mathfrak{sl}_2)$
coproduct realized in $\calf\otimes \bar\calf$ we set
\begin{equation}
\bar M^{-1}_+M_+=\left(
\begin{array}{cc}
\Delta(\uk^{-\frac{1}{2}}) & (
q^{-1}-q)\Delta(\uf)\Delta(\uk^{1/2})\\
\\
0 & \Delta(\uk^{1/2})
\end{array}
\right)\label{2.2}
\end{equation}  
\begin{equation}
\bar M_-^{-1}M_-=\left(
\begin{array}{cc}
\Delta(\uk^{1/2}) & 0\\
\\
(q-q^{-1})\Delta(\uk^{-1/2})\Delta(\ue) & \Delta(\uk^{-1/2})
\end{array}
\right)
\label{2.3}
\end{equation}
 where $M_\pm(\bar M_\pm)$ are the Gauss components of $M(\bar M)$
defined by
\[
q^{3/2}M=M_+M^{-1}_-,\ \ q^{3/2}\bar M^{-1}=\bar M^{-1}_+\bar
M_-,
\]
with the same diagonal elements in $M_+$ and $M_-^{-1}$ (as well
as in $\bar M_+$ and $\bar M_-^{-1}$).  
Noting that $M^{\pm 1}_\pm$ and $\bar M^{\mp 1}_\pm$ satisfy
identical exchange relations (\cite{FHT2}, \cite{FHT3}), we
parametrize them in the same way in terms of generators $X$ and
$\bar X$ of the corresponding QUEA as the products
(\ref{2.2}),(\ref{2.3})
are expressed in terms of $\Delta(\ux)$. As a result we obtain
\[
\Delta(\uk^{\pm 1/2})=q^{\pm 1/2(H+\bar H)}, \ \Delta (\ue)=E+q^H\bar
E,\ \ \Delta(\uf)=F q^{-\bar H}+\bar F
\]
(where the two copies of $U_q(\mathfrak{sl}_2)$ labelled by $X$ and $\bar
Y$ commute~: $[X,\bar Y]=0)$. The second and the third QUEA are
generated by $A,A'$ and $B,B'$  ($[A^{(\prime)},B^{(\prime)}]=0$)
where 
\[
A=a^2_\alpha\bar a^\alpha_2,\ \ A'=a^1_\alpha \bar a^\alpha_1
\Rightarrow [A,A']=[p+\bar p],\ \  q^{p+\bar p} A=Aq^{p+\bar p-2}
\]
\[
B=a^1_\alpha\bar a^\alpha_2,\ \ B'=-a^2_\alpha \bar
a^\alpha_1\Rightarrow [B,B']=[p-\bar p],\ \ q^{p-\bar p}
B=Bq^{p-\bar p+2}
\]
We shall denote the QUEA generated by $\Delta(\ux)$ and by $B$
and $B'$ by $U_q(\mathfrak{sl}_2)_\Delta$ and $U_q(\mathfrak{sl}_2)_B$,
respectively. In order to prove the $U_q(\mathfrak{sl}_2)_\Delta$
invariance of $a^i_\alpha\bar a^\alpha_j$ we use the exchange
relations
\[
q^H a^i_\alpha=a^i_\alpha q^{H+\delta^1_\alpha-\delta^2_\alpha},\
q^{\bar H}\bar a^\alpha_i=\bar a^\alpha_iq^{\bar
H-\delta^1_\alpha+\delta^2_\alpha},
\]
\[
[E,a^i_\alpha]=\delta^2_\alpha a^i_1q^H,\ q^{3-2\alpha}
\bar E\bar
a^\alpha_i-\bar a^\alpha_i \bar E=-\delta^\alpha_1\bar a^2_i,\
\alpha=1,2,
\]
\[
F a^i_\alpha-q^{2\alpha-3}a^i_\alpha F=\delta^1_\alpha a^i_2,\ \
[\bar a^\alpha_i,\bar F]=\delta^\alpha_2q^{-\bar H}\bar a^1_i
\]
which imply $[\Delta(\ux),a^i_\alpha\bar a^\alpha_j]=0$
   for $\ux=\ue,\uf,\uk, i,j=1,2$.
The results of \cite{D-VT} (propositions 1,2,4) can be summarized
as follows.
\addtocounter{theorem}{-1}
\begin{theorem}
$(a)$ The set of vectors in $\calh$ invariant under the pair of
mutually  commuting QUEA $U_q(\mathfrak{sl}_2)_\Delta$ and
$U_q(\mathfrak{sl}_2)_B$ spans a $2h-1$ dimensional space $\calh_I$ with
basis $\{\vert n+1>_I=(A')^{[n]}\vert 1,0>\otimes \vert 1,0>$,
$n=0,1,\dots,2h-2\}(\subset\calh_I)$ where (for
$A'_\alpha=a^1_\alpha\bar a^\alpha_1$,  $\alpha=1,2$ (no
summation in $\alpha$))
\[
(A')^{[n]}=\frac{1}{[n]!}(A')^n=\sum^{n-m}_{\ell=m}q^{\ell(n-\ell)}
(A'_1)^{[\ell]}(A'_2)^{[n-\ell]},m=\max(0,n-h+1).
\]
$(b)$ The operators $A_1=a^2_1\bar a^1_2$ and $A_2=a^2_2\bar
a^2_2$, as well as $A'_\alpha$ satisfy $A^{(\prime)}_2
A^{(\prime)}_1=q^2 A^{(\prime)}_1
A^{(\prime)}_2$ and $(A^{(\prime)}_\alpha)^h=0$ in
$\calh=\calf\otimes \bar \calf$ implying $(A^{(\prime)})^h=0$ 
for $A^{(\prime)}$ standing for $A$ or $A'$; furthermore, the
basis ($\vert n>_I$) is characterized by
\[
A\vert n>_I=[n]\vert n -1>_I,\ \ ([p]-[n])\vert n>_I=0.
\]
$(c)$ Each of the generalized homologies of the nilpotent
operator $A$ in $\calh_I$ is one-dimensional and given by
\[
H_{(n)}(\calh_I,A)\simeq \{\mathbb C\vert
n>\},\ n=1,2,\dots,h-1.\]
\end{theorem}
When $q$ is a root of the unity, $U_q(\mathfrak{sl}_2)$ has a huge Hopf
ideal for which the quotient $\tilde U_q(\mathfrak{sl}_2)$ is a finite
dimensional Hopf algebra (the {\sl reduced} QUEA). Here we have
$q^h=-1$ and it is not hard to see that the actions of
$U_q(\mathfrak{sl}_2)_\Delta$ and $U_q(\mathfrak{sl}_2)_B$ on $\calh$ are in
fact actions of the corresponding finite dimensional quotients.
In the following, $\calu_q$ will denote their tensor product. It is this finite-dimensional Hopf algebra $\calu_q$
which acts on $\calh$ and $\calh_I$ is the subspace of
$\calu_q$-invariant vectors.

\Section{Necessity to enlarge $\calh$}

In short, one has a vector space $\calh$ on which act a Hopf
algebra $\calu_q$ and a nilpotent endomorphism $A$
satisfying $A^h=0$. The action of the algebra $\calu_q$ commutes
with $A$, i.e. one has on $\calh$ : $ [A,X]=0,\ \ \ \forall 
X\in \calu_q$.
It follows that the subspace $\calh_I$ of $\calu_q$-invariant
vectors in $\calh$ is stable by $A$, i.e. $A(\calh_I)\subset
\calh_I$. Thus $(\calh_I,A)$ is an $h$-differential subspace of
the $h$-differential vector space $(\calh,A)$ and it turns out
that the ``interesting object" (the physical space) is the
generalized homology of $(\calh_I,A)$. We would like to avoid the
restriction to the invariant subspace $\calh_I$ that is, in
complete analogy with the BRS methods, we would like to define an
extended $h$-differential space in such a way that the
$\calu_q$-invariance is captured by its $h$-differential in the
sense that it
has the same generalized homology as $(\calh_I,A)$.\\
The most natural thing to do is to try to construct a nilpotent
endomorphism $Q$ of $\calh$ with $Q^h=0$ such that its
generalized homology coincides with the one of $A$ on $\calh_I$
i.e. such that one has
  \begin{equation}
H_{(n)}(\calh,Q)=H_{(n)}(\calh_I,A),\ \ \ \forall n\in
\{1,\dots,h-1\}
\label{3.2}
\end{equation}
   where $H_{(n)}(\calh,Q)=\ker(Q^n)/\im(Q^{h-n})$.
Unfortunately this is  not possible. Indeed let $Q$ be
a nilpotent endomorphism of $\calh$ as above and let us decompose
$\calh$ into irreducible factors\footnote{We are using here the
terminology of \cite{Gr}, (Chapter XV, \S 3); (maximal) 
indecomposable factors would perhaps be better than irreducible
factors.} for $Q$ \cite{Gr}. One obtains an
isomorphism\linebreak[4] $\calh\simeq \displaystyle{\oplusinf^h_{n=1}} 
\mathbb C^n\otimes \mathbb C^{m_n}$, $Q\simeq \displaystyle{\oplusinf^h_{n=1}}
Q_n\otimes \id_{\mathbb C^{m_n}}$ with
\[
Q_n=\left(
\begin{array}{ccccccc}
0 & 1 & 0 & . & . & . & 0\\
. & . & . & . &  & &.\\
. &   & . & . & . &  & .\\
. &   &   & . & . & . & .\\
. &  &   &   & .  & . & 0\\
. & & & & &.& 1\\
0 & . &. &. & . & . & 0\\
\end{array}
\right)\in M_n(\mathbb C)\ \ \ \  ,(Q_1=0),
\]
where $m_n$ is the multiplicity of the irreducible
representation in $\mathbb C^n$. The above decomposition is also
the Jordan normal form of $Q$. One has
$\dim(\calh)=\sum^h_{n=1}nm_n$ and it is easy to compute
$\dim(H_{(n)}(\calh,Q))$ in terms of the multiplicities $m_n$.
The result is (Proposition 2 of \cite{D-V2}) 
  \begin{equation}
\dim(H_{(n)}(\calh,Q))=\sum^n_{j=1}\sum^{h-j}_{i=j}m_i=
\dim(H_{(h-n)}(\calh,Q))
\label{3.3}
\end{equation}
for $1\leq n\leq h/2$.
On the other hand we know from \cite{D-VT} that one has
$\dim(H_{(n)}(\calh_I,A))=1$ for $1\leq n\leq h-1$. This
implies, by using (\ref{3.2}) and (\ref{3.3}), that one must have either
$m_{h-1}=1$ and $m_n=0$ for $1\leq n\leq h-2$ or $m_1=1$ and
$m_n=0$ for $2\leq n\leq h-1$. It follows that one must have
either $\dim(\calh)=hm_h+h-1$ or $\dim(\calh)=hm_h+1$. However
$\dim(\calh)=\dim(\calf\otimes \bar\calf)=h^4$ as in Section 2
is  not compatible with the above estimates. The
same conclusion would hold for other natural choices for
$\calh$. It is worth noticing here that $\calh_I$ is perfectly
compatible with the first possibility since
$\dim(\calh_I)=2h-1=h+h-1$.\\

\Section{A minimal canonical construction}

We recall that $A$ is a nilpotent endomorphism of $\calh$ with
$A^h=0$ and $A(\calh_I)\subset \calh_I$. Let us
define the graded vector space
$\calh^\bullet=\displaystyle{\oplusinf_{n\geq 0}}\calh^n$ by
$\calh^0=\calh$, $\calh^n=\calh/\calh_I$ for $1\leq n\leq h-1$
and $\calh^n=0$ for $n\geq h$. One then defines an endomorphism
$d$ of degree 1 of $\calh^\bullet$ by setting
$d=\pi:\calh^0\rightarrow \calh^1$ where $\pi:\calh\rightarrow
\calh/\calh_I$ is the canonical projection,
$d=\id:\calh^n\rightarrow \calh^{n+1}$ for $1\leq n\leq h-2$
where $\id$ is the identity mapping of $\calh/\calh_I$ onto
itself and $d=0$ on $\calh^n$ for $n\geq h-1$. One has $d^h=0$
and therefore $(\calh^\bullet,d)$ is an $h$-complex, so its generalized (co)homology is
graded  $H_{(k)}(\calh^\bullet,d)=
\displaystyle{\oplusinf_{n\geq 0}} H^n_{(k)}(\calh^\bullet,d)$ 
with
\[
H^n_{(k)}(\calh^\bullet,d)=\ker(d^k:\calh^n\rightarrow
\calh^{n+k})/d^{h-k}(\calh^{n+k-h})
\]
It is given by the
following proposition.
\begin{proposition}
One has $H^n_{(k)}(\calh^\bullet,d)=0$ for $n\geq 1$
and\linebreak[4]
$H^0_{(k)}(\calh^\bullet,d)=\calh_I$, $\forall k\in \{1,\dots,
h-1\}$.
\end{proposition}
\noindent \underbar{Proof.} Let $\varphi\in \calh^n$ for $n\geq
1$ be such that $d^k\varphi=0$. Then either $\varphi=0$ or
$k\geq  h-n$. In the latter case one has
$\varphi=d^n\psi=d^{h-k}(d^{n+k-h}\psi$) which implies that the
class of $\varphi$ vanishes in $H^n_{(k)}(\calh^\bullet,d)$.
This proves $H^n_{(k)}(\calh^\bullet,d)=0$ for $n\geq 1$. Let
$\psi\in \calh^0=\calh$ be such that $d^k\psi=0$. Then, by
definition this is equivalent to $\pi(\psi)=0$ i.e. $\psi\in\calh_I$
which achieves the proof of the proposition. $\square$\\
It is worth noticing here that given the vector space $\calh$
together with the subspace $\calh_I$, the $h$-complex
$(\calh^\bullet,d)$ is characterized (uniquely up to an
isomorphism) by the following universal property (the proof of
which is straightforward).
\begin{proposition}
Any linear mapping $\alpha:\calh\rightarrow \calc^0$ of $\calh$
into the subspace $\calc^0$ of elements of degree $0$ of an
$h$-complex $(\calc^\bullet,d)$ which satisfies $d\circ
\alpha(\calh_I)=0$ extends uniquely as a
homomorphism
$\bar\alpha:(\calh^\bullet,d)\rightarrow (\calc^\bullet,d)$ of
$h$-complexes.
\end{proposition}
We now use this universal property to extend $A$ to
$\calh^\bullet$.
\begin{proposition}
The endomorphism $A$ of $\calh=\calh^0$ has a unique extension
to $\calh^\bullet$, again denoted by $A$, as a homogeneous
endomorphism of degree $0$ satisfying $Ad-q^2\ dA=0$. On
$\calh^\bullet$, one has $A^h=0$ and $(d+A)^h=0$.
\end{proposition}
\noindent \underbar{Proof.} Since $A(\calh_I)\subset \calh_I$,
one can apply the universal property (Proposition 2) for
$\alpha=A:\calh\rightarrow \calh^0$ and one obtains a unique
homomorphism $\bar A:\calh^\bullet\rightarrow \calh^\bullet$ of
$h$-complexes  extending $A$. One has $\bar Ad=d\bar A$ which is
equivalent to $Ad-q^2dA=0$ for $A=q^{2D}\bar A=\bar A q^{2D}$
where $D$ denotes the degree in $\calh^\bullet$. Again by
uniqueness in Proposition 2, one has $\bar A^h=0$ which is
equivalent to $A^h=0$ on $\calh^\bullet$. It follows
from $Ad-q^2dA=0$ and from the fact that $q^2$ is a primitive
$h$-root of the unity that one has $(d+A)^h=d^h+A^h$ which
implies $(d+A)^h=0$. $\square$\\
Thus $Q= d+A$ is an $h$-differential.
The main result of this section, Theorem 1, states that the
generalized homology $H_{(k)}(\calh^\bullet, Q)$ coincides
with $H_{(k)}(\calh_I,A)$. In order to prove the result we shall
need the following construction and lemma, (Lemma 1). Let
$\cale$ be a vector space equipped with a nilpotent endomorphism
$L$ satisfying $L^h=0$, (i.e. $(\cale,L)$ is an $h$-differential
space), and let
$\cale^\bullet=\displaystyle{\oplusinf_n}\cale^n$ be the graded
vector space defined by setting $\cale^n=\cale$ for $0\leq n\leq
h-1$ and $\cale^n=0$ otherwise. Let $\delta$ and $\call$ be the
 endomorphisms of
the vector space $\cale^\bullet$ defined by setting
$\delta(\psi)_n=\psi_{n-1}$ and  $\call(\psi)_n=q^{2n}L(\psi_n)$
for $0\leq n\leq h-1$
with
$\psi=\displaystyle{\oplusinf_n}\psi_n$ ($\psi_n\in\cale^n$).
One has $(\delta+\call)^h=0$ because 
$\call\delta-q^2\delta\call=0$, $\delta^h=\call^h=0$. Thus
$(\cale^\bullet,\delta+\call)$ is an $h$-differential vector space.
\begin{lemma}
One has $H_{(k)}(\cale^\bullet,\delta+\call)=0$ for $1\leq k\leq h-1$.
\end{lemma}
\noindent \underbar{Proof.} In view of Lemma 3 in \cite{D-V2} it
is sufficient to prove that one has\linebreak[4] 
$H_{(1)}(\cale^\bullet,\delta+\call)=0$. So let $\psi$ be such that
$(\delta+\call)(\psi)=0$. By definition this means $L(\psi_0)=0$ and
$\psi_{n-1}+q^{2n}L(\psi_n)=0$ for $1\leq n\leq h-1$ which is
equivalent to $\psi_k=(-1)^{h-1-k}\ q^{h(h-1)-k(k+1)}\
L^{h-1-k}(\psi_{h-1})$ for $0\leq k\leq h-1$ (since
$L(\psi_0)=0$ follows then from $L^h=0$). On the other hand let
$\varphi=\displaystyle{\oplusinf_n}\varphi_n\in \cale^\bullet$
be defined by $\varphi_0=(-1)^{h-1}q^{h(h-1)}\psi_{h-1}=
\psi_{h-1}$ and
 $\varphi_n=0$ otherwise; then one has $(\delta+\call)^{h-1}
 (\varphi)=\psi$.
This proves that $H_{(1)}(\cale^\bullet,\delta+\call)=0$ and implies the
result. $\square$
\begin{theorem}
The generalized $Q$-homology of $\calh^\bullet$ coincides with
the generalized $A$-homology of $\calh_I$, i.e. one has 
$H_{(k)}(\calh^\bullet,Q)=H_{(k)}(\calh_I,A)$
for $1\leq k\leq h-1$.
\end{theorem}
\noindent\underbar{Proof.} Let us consider the previous
$h$-differential vector space $(\cale^\bullet,\delta+\call)$ for the
choices $\cale=\calh/\calh_I$ and $L=A_\pi$. One defines a
surjective linear mapping $\beta$ of $\calh^\bullet$ onto 
$\cale^\bullet$ by
setting $\beta=\id : \calh^n\rightarrow \cale^n$ for $1\leq
n\leq h-1$ and $\beta=\pi:\calh^0\rightarrow \cale^0$ where
$\id$ is the identity mapping of $\calh/\calh_I$ onto itself and
where $\pi$ is the canonical projection of $\calh$ onto
$\calh/\calh_I$. The kernel of $\beta$ is obviously $\calh_I$ so
one has a short exact sequence
\[
0\rightarrow \calh_I\stackrel{\alpha}{\rightarrow}
\calh^\bullet\stackrel{\beta}{\rightarrow} \cale^\bullet\rightarrow 0
\]
where $\alpha$ is the composition of inclusions $\calh_I\subset
\calh=\calh^0\subset \calh^\bullet$. It is straightforward to
verify that one has $\alpha \circ A=(d+A)\circ \alpha$ and
$\beta\circ(d+A)=(\delta+\call)\circ\beta$ so one has in fact a short exact
sequence of $h$-differential vector spaces
\[
0\rightarrow
(\calh_I,A)\stackrel{\alpha}{\rightarrow}(\calh^\bullet,Q)
\stackrel{\beta}{\rightarrow}(\cale^\bullet,\delta+\call)\rightarrow 0
\]
By using  Lemma 1 above and Lemma 2 of \cite{D-V2}, one obtains
the exact sequences
\[
0\stackrel{\partial}{\rightarrow}H_{(k)}(\calh_I,A)
\stackrel{\alpha_\ast}{\rightarrow}
H_{(k)}(\calh^\bullet,Q)\stackrel{\beta_\ast}{\rightarrow} 0
\]
where $\alpha_\ast$ and $\beta_\ast$ are induced by $\alpha$ and
$\beta$ and where $\partial$ is the connecting homomorphism
\cite{D-V2}, \cite{KW}. Thus $\alpha_\ast$ is an isomorphism 
which 
allows the canonical identifications of Theorem 1. $\square$\\
\noindent\underbar{Remark 1.} Notice that the content of this
section does only depend on the data $(\calh,A,\calh_I)$ where
$(\calh,A)$ is an $h$-differential vector space and $\calh_I$ is a
subspace of $\calh$ invariant by $A$. The same remark applies to
Section 3, except that, of course, the specific dimensions are
also involved there.

\Section{Extension to  Hochschild cochains}

Although the construction of last section is quite optimal, it
is lacking \linebreak[4] a ``geometrico-physical" interpretation. 
Our aim in
the following is to cure that by casting the construction in a
form which is closer to the BRS formulation in gauge theory or
in constrained systems. To this end we recall that $\calh$ is a
representation space for the Hopf algebra
$\calu_q$ and that $\calh_I$
is the subspace of $\calu_q$-{\sl invariant elements} of
$\calh$. An element $\Psi\in \calh$ is said to be $\calu_q$-{\sl
invariant}, or simply {\sl invariant} when no confusion arises,
if it satisfies
  \begin{equation}
X\Psi=\Psi \varepsilon (X),\ \ \ \forall X\in \calu_q
\label{5.1}
\end{equation}
   where $\varepsilon$ denotes the counit of $\calu_q$. It turns
out that (\ref{5.1}) has a natural interpretation in terms of Hochschild
cohomology. To see this, we equip $\calh$ with a structure of
bimodule over $\calu_q$. One already has a structure of left
$\calu_q$-module on $\calh$ given by the representation of
$\calu_q$ in $\calh$. We equip $\calh$ with a structure of right
$\calu_q$-module by using the scalar representation given by the
counit $\varepsilon$. Since one obviously has
$(X\Psi)\varepsilon (Y)=X(\Psi\varepsilon(Y))$ for any
$\Psi\in \calh$ and $X,Y\in \calu_q$, $\calh$ is a bimodule. One
can introduce the graded space
$C(\calu_q,\calh)=\displaystyle{\oplusinf_{n\geq 0}}
C^n(\calu_q,\calh)$ of $\calh$-valued Hochschild cochains of
$\calu_q$, where $C^n(\calu_q,\calh)$ is the vector space of all
linear mappings of $\stackrel{n}{\otimes}\calu_q$ into $\calh$,
(i.e. $n$-linear mappings of $(\calu_q)^n$ into $\calh$).
Equipped with the Hochschild differential $d_H$,
$C(\calu_q,\calh)$ is a complex and the $\calh$-valued
Hochschild cohomology of $\calu_q$, 
$H(\calu_q,\calh)=\displaystyle{\oplusinf_{n\geq 0}}
H^n(\calu_q,\calh)$,  is the homology of this complex. Now the
condition (\ref{5.1}) for $\Psi$ to be in $\calh_I$ also reads
$d_H\Psi=0$. Therefore $\calh_I$ identifies with
$H^0(\calu_q,\calh)$.
However, except for $h=2$, one cannot mix reasonably the
Hochshild differential $d_H$ satisfying $d^2_H=0$ with (an
extension of) the nilpotent $A$ satisfying $A^h=0$. Fortunately,
there is an $h$-differential $d$ on $C(\calu_q,\calh)$ which
coincides with $d_H$ in degree 0. This $d$ was introduced in
\cite{D-VK} (with the notation $d=d_{q^2}$) and was analysed in details
in \cite{D-V2} (with the notation $d=d_1$; see Remark 3 below). 
It is given for $\omega\in C^n(\calu_q,\calh)$
by
  \begin{eqnarray}
d(\omega)(X_0,\dots,X_n) & = & X_0\omega(X_1,\dots,X_n)\nonumber\\
& + & \sum^n_{k=1}
q^{2k}\omega(X_0,\dots,(X_{k-1}X_k),\dots,X_n)\nonumber\\
& - & q^{2n}\omega(X_0,
\dots, X_{n-1})\varepsilon(X_n).
\label{5.2}
\end{eqnarray}
 \begin{lemma}
Let $\Psi\in \calh=C^0(\calu_q,\calh)$; the following
conditions $(i)$, $(ii)$ and $(iii)$ are equivalent\\
$(i)$ $d^k(\Psi)=0$ for some $k$ with $1\leq k\leq h-1$\\
$(ii)$ $\Psi\in\calh_I$\\
$(iii)$ $d^n(\Psi)=0$ for any $n\in\{1,\dots,h-1\}$.
\end{lemma}
\noindent\underbar{Proof.} We know that $\Psi\in \calh_I$ is
equivalent to $d_H\Psi=0$ and, since $d=d_H$ on
$C^0(\calu_q,\calh)$ this is equivalent to $d\Psi=0$. The
implication ($ii$)$\Rightarrow$($iii$) follows (since
$d\Psi=0\Rightarrow d^n\Psi=0$ for $n\geq 1$). The
implication ($iii$)$\Rightarrow$($i$) is clear. It remains to
show the implication ($i$)$\Rightarrow$($ii$) to achieve the
proof of the lemma. By induction on $n$ and by using 
definition (\ref{5.2}) one has for $\Psi\in C^0(\calu_q,\calh)$
  \begin{equation}
d^n\Psi(\bbbone,\dots,\bbbone,X)=(1+q^2)\dots
(1+q^2+\dots+q^{2(n-1)})d\Psi(X)
\label{5.3}
\end{equation}
   for any $n\geq 1$, $X\in \calu_q$ where $\bbbone$ is the unit of
$\calu_q$. Formula (\ref{5.3}) shows that $d^k\Psi=0$ for some
$k\in\{1,\dots,h-1\}$ implies $d\Psi=0$, i.e. $\Psi\in \calh_I$
and thus the implication ($i$)$\Rightarrow$ ($ii$). $\square$\\
As an easy consequence of this lemma one obtains the following
result.
\begin{proposition}
The $h$-complex $(\calh^\bullet,d)$ can be canonically
identified with the $h$-subcomplex of $(C(\calu_q,\calh),d)$
generated by $\calh$ (that is with $\calh\oplus d\calh\oplus\dots\oplus d^{h-1}\calh\subset
C(\calu_q,\calh)$ for the $h$-differential $d$).
\end{proposition}
\noindent\underbar{Proof.} In view of Proposition 2 (i.e. in
view of the universal property of $(\calh^\bullet,d)$), the
identity mapping of $\calh$ onto itself extends uniquely as a
homomorphism of $h$-complexes of $(\calh^\bullet,d)$ into
$(C(\calu_q,\calh),d)$. Lemma 2 then implies that this
homomorphism is injective. $\square$\\
Thus one has $\calh^\bullet \subset C(\calu_q,\calh)$ and the
$h$-differential $d$ of $C(\calu_q,\calh)$ extends the one of
$\calh^\bullet$; we now extend $A$ to $C(\calu_q,\calh)$. 
\begin{proposition}
Let us extend $A$ to $C(\calu_q,\calh)$ as a homogeneous
endomorphism $\omega\mapsto (A\omega)$ of degree $0$ by setting
\[
(A\omega)(X_1,\dots,X_n)=q^{2n}A\omega(X_1,\dots,X_n)
\]
for $\omega\in C^n(\calu_q,\calh)$ and $X_i\in\calu_q$. On 
$C(\calu_q,\calh)$one has
 $Ad-q^2dA=0$, $A^h=0$ and $(d+A)^h=0$.
\end{proposition}
\noindent\underbar{Proof.} Consider first the extension
$\bar A$ defined by $(\bar
A\omega)(X_1,\dots,X_n)=A\omega(X_1,\dots,X_n)$ for $\omega\in
C^n(\calu_q,\calh)$ and $X_i\in\calu_q$. Then, by using the
fact that the action of $\calu_q$ on $\calh$ commutes with $A$
(see Section 2), one obtains $\bar Ad=d\bar A$; more
generally, if the {\sl cofaces} of $C(\calu_q,\calh)$
\cite{D-V2}
\[
f_\alpha:C^n(\calu_q,\calh)\rightarrow
C^{n+1}(\calu_q,\calh),\ \  \alpha\in\{0,\dots,n+1\}
\]
are defined by
\[
\begin{array}{lll}
(f_0\omega)(X_0,\dots,X_n)& = & X_0\omega(X_1,\dots,X_n),\\
\\
(f_i\omega)(X_0,\dots,X_n)& = & \omega(X_0,\dots,(X_{i-1}X_i),
\dots,X_n)
\end{array}
\]
for  $i\in\{1.\dots,n\}$ and
\[
(f_{n+1}\omega)(X_0,\dots,X_n) = \omega(X_0,\dots,X_{n-1})
\varepsilon (X_n),
\]
then one has $
\bar Af_\alpha  = f_\alpha\bar A,
\forall \alpha \in \{0,\dots,n+1\}.$
This implies that $A=q^{2D}
\bar A=\bar A q^{2D}$ satisfies in particular $Ad-q^2dA=0$,
$D$ being the cochain's degree; more generally $A$ satisfies
$Af_\alpha=q^2f_\alpha A$. Furthermore $\bar A^h=0$, which
is straightforward, implies $A^h=0$. The last equality
$(d+A)^h=0$ follows by the same argument as in the proof of
Proposition 3. $\square$\\
We have now extended to $C(\calu_q,\calh)$ the whole structure
defined on $\calh^\bullet$ in the previous section. Indeed the uniqueness in Proposition 3 implies
that $A$ defined on $C(\calu_q,\calh)$ in Proposition 5 is an
extension of $A$ defined on $\calh^\bullet$ in Proposition 3.
One then extends to $C(\calu_q,\calh)$ the definition of $Q$ by
setting again $Q=d+A$.
Next section will be devoted to the formulation and the
discussion of the appropriate extension to $C(\calu_q,\calh)$ of
Theorem 1.\\
\noindent\underbar{Remark 2}. Lemma 2 implies : 
$H^0_{(k)}(C(\calu_q,\calh),d)=H^0(\calu_q,\calh)$,\linebreak[4] $\forall
k\in\{1,\dots,h-1\}$. This is a special case of Theorem 4
(1) of \cite{D-V2} which reads here~:
$H^{hr}_{(k)}(C(\calu_q,\calh),d)=H^{2r}(\calu_q,\calh)$, 
$H^{h(r+1)-k}_{(k)}(C(\calu_q,\calh),d)=H^{2r+1}
(\calu_q,\calh)$,
and $H^n_{(k)}(C(\calu_q,\calh),d)=0$ otherwise.

\noindent\underbar{Remark 3}. Given primitive $h$-th root of
the unity $q^2$, one can construct on $C(\calu_q,\calh)$ several
$h$-differentials of degree 1 which coincide with the Hochschild
differential when $q^2=-1$. A whole sequence $(d_n)_{n\in \mathbb
N}$ of such $h$-differentials has been introduced  in \cite{D-V2}
where their generalized cohomologies were computed in terms of
the ordinary Hochschild cohomology. For the case of $d_0$, a
$h$-differential which has be considered by several authors
(\cite{Kap}, \cite{KW}), this computation has been done
independently by Kassel and Wambst by using very interesting
generalizations of concepts of homological algebra. However here
only $d_1=d$ should be used for $h>2$ because $d_0$ does not
coincide with the Hochschild differential $d_H$ on $\calh$
and, on the other hand, although $d_n$ coincides with $d_H$ on
$\calh$ for $n\geq 1$ it also coincides with $d_H$ on the
1-cochains for $n\geq 2$ and thus $d^2_n$ vanishes on $\calh$
whenever $n\geq 2$.

\Section{Generalized homology of $Q$ on $C(\calu_q,\calh)$}

As explained in \cite{D-V2} (see Remark 2 above) the spaces
$H^n_{(k)}(C(\calu_q,\calh),d)$ can be computed in terms of the
Hochschild cohomology $H(\calu_q,\calh)$. In particular, one sees
that $H^n_{(k)}(C(\calu_q,\calh),d)$ does not generally vanish
for $n\geq 1$. This implies that one cannot expect for the
generalized homology of $Q$ on $C(\calu_q,\calh)$ such a simple
result as the one given by Theorem 1 for the generalized homology
of $Q$ on $\calh^\bullet$. Nevertheless, in view of Lemma 2, one
has
$H^0_{(k)}(C(\calu_q,\calh),d)=\calh_I=H^0_{(k)}(\calh^\bullet,d)$
and therefore one may expect 
$H^0_{(k)}(C(\calu_q,\calh),Q)=H_{(k)}(\calh_I,A)(=H_{(k)}
(\calh^\bullet,Q))$. In fact, this is essentially true. However
some care must be taken because $Q$ is not homogeneous so
$H_{(k)}(C(\calu_q,\calh),Q)$ is not a graded vector space.
Instead of a graduation, one has an increasing filtration
$F^nH_{(k)}(C(\calu_q,\calh),Q)$, ($n\in \mathbb Z$),
with
$F^nH_{(k)}(C(\calu_q,\calh),Q)=0$ for $n<0$ and where, for 
$n\geq
0$, $F^nH_{(k)}(C(\calu_q,\calh),Q)$ is the canonical image  in
$H_{(k)}(C(\calu_q,\calh),Q)$ of $\ker(Q^k)\cap
\displaystyle{\oplusinf^{r=n}_{r=0}}C^r(\calu_q,\calh)$. There is
an associated graded vector space
\[
^{\mathrm gr}H_{(k)}(C(\calu_q,\calh),Q)=\oplusinf_n
F^nH_{(k)}(C(\calu_q,\calh),Q)/F^{n-1}H_{(k)}(C(\calu_q,\calh),Q)
\]
which here is $\mathbb N$-graded. One has
\[
F^0H_{(k)}(C(\calu_q,\calh),Q)=^{\mathrm
gr}H_{(k)}^0(C(\calu_q,\calh),Q)
\]
and it is this space which is the correct version of the
$H^0_{(k)}(C(\calu_q,\calh),Q)$ above in order to identify
$H_{(k)}(\calh_I,A)$ in the generalized homology of $Q$ on
$C(\calu_q,\calh)$.
\begin{theorem}
The inclusion $\calh^\bullet\subset C(\calu_q,\calh)$ induces the
 isomorphisms 
\[
H_{(k)}(\calh^\bullet,Q)\simeq
F^0H_{(k)}(C(\calu_q,\calh),Q)\  \mathrm{for}\   1\leq k\leq h-1.
\]
In
particular, with obvious identifications, one has
\[
F^0H_{(k)}(C(\calu_q,\calh),Q)=H_{(k)}(\calh_I,A),\ \ \ \forall
k\in \{1,\dots,h-1\}.
\]
\end{theorem}
\noindent\underbar{Proof.} 
Let $\Psi\in \calh$ be such that one has $Q^k\Psi=0$ for some
$k$ ($1\leq k\leq h-1$). By expanding $(d+A)^k\Psi=0$, one obtains $d^k\Psi=0$ for the highest degree and
$A^k\Psi=0$ for the lowest degree. In view of Lemma 2 this is
equivalent to $\Psi\in\calh_I$ and $A^k\Psi=0$; this conversely
implies $Q^k\Psi=0$. Thus $Q^k\Psi=0$ for $\Psi\in \calh$ is
equivalent to $\Psi\in\calh_I$ and $A^k\Psi=0$. On the other
hand $\Psi=Q^{h-k}\Phi$ for $\Psi\in \calh$ implies $\Phi\in
\calh$ and $d^{h-k}\Phi=0$ which by using again Lemma 2 implies
$\Phi\in\calh_I$ and $\Psi=A^{h-k}\Phi=Q^{h-k}\Phi$\linebreak[4]
$(\in\calh_I)$.
This means that one has canonically:
\[
F^0H_{(k)}(C(\calu_q,\calh),Q)=H_{(k)}(\calh_I,A)=H_{(k)}(\calh^\bullet,Q)
\]
which completes the proof of Theorem 2. $\square$\\
If one compares this construction involving Hochschild cochains 
with the construction of Section 4, what has
been gained here besides the explicit occurrence of the quantum
gauge aspect is that the extended space $C(\calu_q,\calh)$ is a
tensor product $\calh\otimes\calh'$ of the original space
$\calh$ with the tensor algebra $\calh'=T(\calu^\ast_q)$ of the
dual space of $\calu_q$. The factor $\calh'$ can thus be
interpreted as the state space for some generalized ghost.
What has been lost is the minimality of the generalized
homology, i.e. besides the ``physical" $H_{(k)}(\calh_I,A)$, the
generalized homology of $Q$ on $\calh\otimes \calh'$ contains
some other non trivial subspace in contrast to what happens on
$\calh^\bullet$. In the usual homological (BRS) methods however
such a ``non minimality" also occurs. Indeed, for instance, in
the homological approach to constrained classical systems, the
relevant homology contains besides the functions on the reduced
phase space the whole cohomology of longitudinal forms
\cite{MD-V}. The same is true for the BRS cohomology of gauge
theory \cite{BRS}, \cite{BC-R}.\\
In the usual situations where one applies the BRS construction
(gauge theory, constrained systems) one has a Lie algebra
$\fracg$ (the Lie algebra of infinitesimal gauge
transformations) acting on some space $\calh$ and what is really
relevant at this stage is the Lie algebra cohomology $H(\fracg,\calh)$ of
$\fracg$ acting on $\calh$. The extended space is then the space
of $\calh$-valued Lie algebra cochains of $\fracg$,
$C(\fracg,\calh)$. This extended space is thus also a tensor
product $\calh\otimes\calh'$ but now $\calh'$ is the exterior
algebra $\calh'=\Lambda\fracg^\ast$ of the dual space of
$\fracg$. That is why this factor can be
interpreted (due to antisymmetry) as a fermionic state space; 
indeed that is the reason why one
gives a fermionic character to the ghost \cite{BRS},
\cite{BC-R}, \cite{Su}, \cite{MD-V}. There is however another way
to proceed in these situations which is closer to what has been
done in our case here. To understand it, we recall that any
representation of $\fracg$ in $\calh$ is also a  
representation of the enveloping algebra $U(\fracg)$ in $\calh$.
Thus $\calh$ is a left $U(\fracg)$-module. Since $U(\fracg)$ is a
Hopf algebra, one can convert $\calh$ into a bimodule for
$U(\fracg)$ by taking as right action the trivial representation
given by the counit. It turns out that the $\calh$-valued Hochschild
cohomology of $U(\fracg)$, $H(U(\fracg),\calh)$, coincides with
the $\calh$-valued Lie algebra cohomology of $\fracg$,
$H(\fracg,\calh)$, i.e. one has \cite{Gui}, \cite{J-LL}:
$H(U(\fracg),\calh)=H(\fracg,\calh)$.
Since it is the latter space which is relevant one can as well
take as extended space the space of $\calh$-valued Hochschild
cochains of $U(\fracg)$, $C(U(\fracg),\calh)$, and then compute
its cohomology. Again this space is a tensor product
$\calh\otimes \calh'$ but now $\calh'=T(U(\fracg)^\ast)$ is a
tensor algebra as in our case.\\
In the above brief discussion of the usual BRS, we
have oversimplified the situation. In general the extended space
contains slightly more than the Lie algebra cochains and the
Chevalley-Eilenberg differential is only a part of the BRS operator.
However the above picture is the essential point. The gist of
our message is the realization that the construction of the
present
paper is in fact very close to the standard BRS procedure. 
The main difference (or
extension) is the occurrence of a nilpotent $Q$ satisfying
$Q^h=0$ with $h>2$ (instead of the usual $Q^2=0$) and,
correspondingly, the occurrence of the generalized homology
(instead of an ordinary homology).

\section*{Acknowledgements}

This work began when I.T. visited l'Institut des Hautes \'Etudes
Scientifiques in Bures-sur-Yvette and le Laboratoire de Physique
Th\'eorique, Universit\'e Paris XI, Orsay (with partial support
from CNRS grant PICS 608) and was continued during his stay at
the Max-Planck Institute for Mathematics in the Sciences in
Leipzig. I.T. thanks these institutions for hospitality and
acknowledges partial support from the Bulgarian National
Foundation for Scientific Research under contract F-827.

\end{document}